\documentclass[a4paper,11pt]{amsart}
\usepackage{amsmath,amsthm,amssymb}
\usepackage[mathscr]{eucal}
 \usepackage{cite}
\usepackage{upgreek}
\usepackage[bookmarks,bookmarksnumbered]{hyperref}
\usepackage{enumerate}
\usepackage{amsmath}
\usepackage{mathrsfs}
\usepackage{blindtext}
\usepackage{scrextend}
\usepackage{enumitem}
\addtokomafont{labelinglabel}{\sffamily}
\usepackage{color}
\usepackage[at]{easylist}

\numberwithin{equation}{section}

\theoremstyle{plain}
\newtheorem{main}{Theorem}
\newtheorem{mcor}[main]{Corollary}

\newtheorem{theorem}{Theorem}[section]

\newtheorem{lemma}[theorem]{Lemma}

\theoremstyle{definition}

\newtheorem*{definition*}{Definition}

\newtheorem{remark}[theorem]{Remark}

\begin{document}

\title[Compact actions whose orbit equivalence relations are not profinite]
{Compact actions whose orbit equivalence relations \\ are not profinite }

\author[A. Ioana]{Adrian Ioana}
\address{Department of Mathematics, University of California San Diego, 9500 Gilman Drive, La Jolla, CA 92093, USA}
\email{aioana@ucsd.edu}

\thanks{The author was supported in part by NSF Career Grant DMS \#1253402.}
\begin{abstract} 
Let $\Gamma\curvearrowright (X,\mu)$ be a measure preserving action of a countable group $\Gamma$ on a standard probability space $(X,\mu)$. 
We prove that if the action $\Gamma\curvearrowright X$ is not profinite and satisfies a certain spectral gap condition, then there does not exist a countable-to-one Borel homomorphism from its orbit equivalence relation to the orbit equivalence relation of any modular action  (i.e., an inverse limit of actions on countable sets).  
As a consequence, we show that if $\Gamma$ is a countable dense subgroup of a compact non-profinite group $G$ such that the left translation action $\Gamma\curvearrowright G$ 
has spectral gap, then $\Gamma\curvearrowright G$ is antimodular and not orbit equivalent to any, {\it not necessarily free}, profinite action.
This provides the first such examples of compact actions,  partially answering a question of Kechris and answering a question of Tsankov.

\end{abstract}

\maketitle

\section{Introduction}

Consider a Borel action $\Gamma\curvearrowright X$ of a countable group $\Gamma$ on a standard Borel space $X$. Then the {\it orbit equivalence relation} $\mathcal R(\Gamma\curvearrowright X)=\{(x_1,x_2)\in X\times X\;|\;\Gamma\cdot x_1=\Gamma\cdot x_2\}$ is a Borel equivalence relation with countable classes. If $\mu$ is a $\Gamma$-invariant Borel probability measure on $X$, then $\mathcal R(\Gamma\curvearrowright X)$ is  a countable probability measure preserving ({\it p.m.p.}) equivalence relation. Conversely,  any countable Borel and p.m.p. equivalence relations can be realized in this way \cite{FM77}.

The study of countable equivalence relations is a central topic in both descriptive set theory and measured group theory (see the surveys \cite{Th06, TS07, Ke18} and \cite{Po07,Fu09,Ga10,Io17}, respectively).  To compare the structure of various equivalence relations one studies class preserving maps. Given equivalence relations $\mathcal R$ and $\mathcal S$ on spaces $X$ and $Y$, a {\it homomorphism} from $\mathcal R$ to $\mathcal S$ is a map $\alpha:X\rightarrow Y$ such that $(\alpha\times\alpha)(\mathcal R)\subset\mathcal S$. 
In the Borel context, one considers Borel homomorphisms, often with some additional properties, such as being one-to-one or a Borel reduction. In the measure theoretic setting, one considers measurable homomorphisms that preserve the underlying measures. In particular, two countable p.m.p. equivalence relations relations $\mathcal R$ and $\mathcal S$ on standard probability spaces $(X,\mu)$ and $(Y,\nu)$ are called isomorphic if there exists an isomorphism of probability spaces $\alpha:X\rightarrow Y$ such that $(\alpha\times\alpha)(\mathcal R)=\mathcal S$,  almost everywhere. Two p.m.p. actions $\Gamma\curvearrowright (X,\mu)$ and $\Delta\curvearrowright (Y,\nu)$ are {\it orbit equivalent} if their orbit equivalence relations are isomorphic.

The aim of this article is to establish a new rigidity phenomenon for countable equivalence relations. Our main result shows that if a p.m.p. action $\Gamma\curvearrowright (X,\mu)$ satisfies a spectral gap condition and its orbit equivalence relation is isomorphic (or admits a countable-to-one homomorphism) to the orbit equivalence relation of some profinite action, then $\Gamma\curvearrowright (X,\mu)$ must be a profinite action.

Before stating our main result in detail, let us first review a few concepts. 
Recall that a p.m.p. action $\Gamma\curvearrowright (X,\mu)$ of a countable group $\Gamma$ on a standard probability space $(X,\mu)$ is called
{\it profinite} if there exists a sequence $\{\mathcal P_n\}$ of finite measurable partitions of $X$ which separates points in $X$, such that each $\mathcal P_n$ is $\Gamma$-invariant. 
 A Borel action $\Delta\curvearrowright Y$ of a countable group $\Delta$  on a standard Borel space $Y$ is called {\it modular} 
if there exists a sequence $\{\mathcal P_n\}$ of countable Borel partitions of $Y$ which separates points in $Y$, such that each $\mathcal P_n$ is $\Delta$-invariant \cite{Hj03}. Note that if $\nu$ is a $\Delta$-invariant ergodic probability measure on $Y$, then after excluding a null set, each partition is finite, the action on each partition is transitive, and the  p.m.p. action $\Delta\curvearrowright (Y,\nu)$ is profinite. Finally, recall that a unitary representation $\pi:\Gamma\rightarrow\mathcal U(\mathcal H)$ on a Hilbert space $\mathcal H$ has {\it spectral gap} if it does not admit almost invariant vectors, i.e., there does not exist a sequence of unit vectors $\xi_n\in\mathcal H$ such that $\|\pi(g)(\xi_n)-\xi_n\|\rightarrow 0$, for all $g\in\Gamma$. A p.m.p. action $\Gamma\curvearrowright (X,\mu)$ has {\it spectral gap} if the Koopman representation  of $\Gamma$ on $L^2_0(X):=L^2(X)\ominus\mathbb C{\bf 1}_X$ does.
The following is our main technical result:

\begin{main}\label{maintech}
Let $\Gamma\curvearrowright (X,\mu)$ be an ergodic p.m.p. action.  
Let $\kappa$ be the Koopman representation of $\Gamma$ on $L^2(X\times X)$ associated to the diagonal action $\Gamma\curvearrowright (X\times X,\mu\times\mu)$. Suppose that the restriction of $\kappa$ to the orthogonal complement of the subspace of $\kappa(\Gamma)$-invariant vectors has spectral gap. 
Assume that there exists a countable-to-one Borel homomorphism from $\mathcal R(\Gamma\curvearrowright X)_{|D}$ to $\mathcal R(\Delta\curvearrowright Y)$, where $D\subset X$ is a measurable set with $\mu(D)=1$ and $\Delta\curvearrowright Y$ is a modular action of a countable group $\Delta$.

Then the action $\Gamma\curvearrowright(X,\mu)$ is profinite. 
\end{main}


Theorem \ref{maintech} implies in particular that if the action $\Gamma\curvearrowright (X,\mu)$ is orbit equivalent to some profinite p.m.p. action $\Delta\curvearrowright (Y,\nu)$, then it must be a profinite action. Thus, roughly speaking, the orbit equivalence relation $\mathcal R(\Gamma\curvearrowright X)$ remembers whether the action it was constructed from is profinite or not. In sharp contrast, if $\Gamma$ is amenable, then  $\mathcal R(\Gamma\curvearrowright X)$ is isomorphic to the unique ergodic hyperfinite p.m.p. equivalence relation, regardless of the action $\Gamma\curvearrowright X$ \cite{Dy58,OW80,CFW81}.

\begin{remark} In the context of Theorem \ref{maintech}, assume that $\Gamma\curvearrowright (X,\mu)$ is orbit equivalent to an (essentially) free profinite p.m.p. action $\Delta\curvearrowright (Y,\nu)$. Suppose additionally that $\Gamma$ has  property (T) or, equivalently (by \cite[Corollary 1.4]{Fu99}), that $\Delta$ has property (T).  Then  by \cite[Theorem A]{Io08}, the action $\Delta\curvearrowright Y$ is orbit equivalence superrigid. Hence, the action $\Gamma\curvearrowright X$ is virtually conjugate to $\Delta\curvearrowright Y$, and so it must be profinite. 
The method of \cite{Io08} relies on a cocycle superrigidity theorem and Zimmer's well-known observation that given any orbit equivalence $\alpha:X\rightarrow Y$ between $\Gamma\curvearrowright X$ and $\Delta\curvearrowright Y$ the formula $\alpha(g\cdot x)=w(g,x)\cdot\alpha(x)$, for $g\in \Gamma$ and $x\in X$,  defines a cocycle $w:\Gamma\times X\rightarrow \Delta$.
On the other hand, if the action $\Delta\curvearrowright Y$ is not free, then one cannot define such a cocycle, and thus the method of \cite{Io08} does not apply.
One of the main novelties of Theorem \ref{maintech} lies in proving that $\Gamma\curvearrowright X$ is profinite {\it without} assuming that the action $\Delta\curvearrowright Y$ is free.
\end{remark}

In the rest of the introduction, we present several consequences of Theorem \ref{maintech}. 
Recall that a p.m.p. action $\Gamma\curvearrowright (X,\mu)$ is called {\it compact}
if the Koopman representation of $\Gamma$ on $L^2(X)$ is a direct sum of finite dimensional representations. 
By a result in \cite{Ma64}, any ergodic compact  action is isomorphic to a left translation action $\Gamma\curvearrowright (G/H,m_{G/H})$ given by $g\cdot (xH)=i(g)xH$, where $G$ is a compact group, $H<G$  is a closed subgroup, $m_{G/H}$ denotes the unique $G$-invariant probability measure on $G/H$, and $i:\Gamma\rightarrow G$ is a homomorphism with dense image. Any profinite p.m.p. action is compact. A left translation action $\Gamma\curvearrowright (G/H,m_{G/H})$ is profinite if and only if $G$ is a profinite group.

By applying Theorem \ref{maintech} to left translation actions on compact groups, we obtain the following:

\begin{mcor}\label{main}
Let $\Gamma$ be a countable dense subgroup of a compact group $G$. Suppose that the left translation action $\Gamma\curvearrowright (G,m_G)$ has spectral gap. Assume that either
\begin{enumerate}
\item The left translation action $\Gamma\curvearrowright (G,m_G)$ is orbit equivalent to a profinite p.m.p. action $\Delta\curvearrowright (Y,\nu)$ of a countable group $\Delta$, or, more generally,
\item There exists a countable-to-one Borel homomorphism from $\mathcal R(\Gamma\curvearrowright G)_{|D}$ to $\mathcal R(\Delta\curvearrowright Y)$, where $D\subset G$ is a measurable set with $m_G(D)=1$ and $\Delta\curvearrowright Y$ is a modular action of a countable group $\Delta$. 
 \end{enumerate}
Then  $G$ is a profinite group.
\end{mcor}

Corollary \ref{main} allows us to answer a question 
originally raised by T. Tsankov in 2009 and revisited by several people since, asking whether every compact action is orbit equivalent to some profinite action. 
 To give some context, recall that in \cite{OP07} 
a notion of  {\it weak compactness} for p.m.p. actions was introduced. This notion, unlike compactness, is an orbit equivalence invariant, and thus lead to a new property for equivalence relations. 
A related property for countable p.m.p. equivalence relations $\mathcal R$ on probability spaces $(X,\mu)$ was defined 
in \cite{Ts09}:   $\mathcal R$  is called {\it profinite} if is the orbit equivalence relation of some profinite action. 
 The main result of \cite{Ts09} provides an intrinsic characterisation of profinite equivalence relations in terms of almost invariant partitions.
  To illustrate the subtle difference between weak compactness and profiniteness, recall that the former requires the existence of a sequence of unit vectors $\xi_n\in L^2(X\times X)$ which are almost invariant under the diagonal action of the full group $[\mathcal R]$  and whose supports converge to the diagonal of $X\times X$ \cite{OP07}. If $\mathcal R$ is the orbit equivalence relation of a profinite action $\Gamma\curvearrowright X$, then one can construct such vectors $\xi_n$ from $\Gamma$-invariant partitions of $X$, and thus $\mathcal R$ is weakly compact.
  On the other hand, it remained open whether every weakly compact equivalence relation (e.g., the orbit equivalence relation of a compact action) is profinite.

We settle the above question negatively here by providing a wide class of compact actions which are not orbit equivalent to any profinite actions. In particular, we show that profiniteness is a strictly stronger property than weak compactness for countable p.m.p. equivalence relations.
Indeed, Corollary \ref{main} implies that if $G$ is compact connected group, then any left translation action $\Gamma\curvearrowright G$ with spectral gap is   not orbit equivalent to any profinite action.
 Following a breakthrough work of J. Bourgain and A. Gamburd \cite{BG06}, the spectral gap property is now known for a large class of translation actions on connected compact groups. Specifically, let $G=SU(d)$, for $d\geq 2$, and $\Gamma<G$ be any dense countable subgroup generated by matrices with algebraic entries. Then the left translation action $\Gamma\curvearrowright G$ has spectral gap \cite{BG06,BG10}. This result has been since generalized in \cite{BdS14} to all connected compact simple Lie groups  $G$.

As a consequence of Corollary \ref{main}, we are also able to partially answer a question raised by A. Kechris in \cite[Section 5 (G)]{Ke05}) asking whether the translation action $\mathbb F_2\curvearrowright G$ arising from an embedding of $\mathbb F_2$ into $G=SO(3)$ is  antimodular.  Recall that a Borel action $\Gamma\curvearrowright X$ is called {\it antimodular} if there does not exist a countable-to-one Borel homomorphism from $\mathcal R(\Gamma\curvearrowright X)$ to $\mathcal R(\Delta\curvearrowright Y)$, for any modular action $\Delta\curvearrowright Y$ \cite{Ke05}. 
Let  $a,b\in\mathbb Z$ such that $0<|a|<b$ and $\frac{a}{b}\not=\pm\frac{1}{2}$. Put $c=b^2-a^2$ and let $\Gamma$ be the subgroup of $G$ generated by the rotations: $$\begin{pmatrix} \frac{a}{b}&-\frac{\sqrt{c}}{b}&0\\\frac{\sqrt{c}}{b}&\frac{a}{b}&0\\0&0&1 \end{pmatrix}\;\;\;\text{and}\;\;\; \begin{pmatrix} 1&0&0\\0&\frac{a}{b}&-\frac{\sqrt{c}}{b}\\0&\frac{\sqrt{c}}{b}&\frac{a}{b}\end{pmatrix}.$$
Then $\Gamma$ is isomorphic to the free group on two generators, $\mathbb F_2$ \cite{Sw94}.
Since the translation action $\Gamma\curvearrowright G$ has spectral gap by \cite{BG06}, Corollary \ref{main} implies that it is antimodular. 
Moreover, one conjecturally expects that the translation action $\Gamma\curvearrowright G$ is antimodular, whenever $\Gamma<G$ is a dense subgroup.
Indeed, a well-known ``spectral gap conjecture" predicts that if the embedding $\Gamma<G$ is dense, then $\Gamma\curvearrowright G$ has spectral gap (see \cite{GJS99}), which by Corollary \ref{main} implies antimodularity. In particular, assuming this spectral gap conjecture, the translation action $\mathbb F_2\curvearrowright G$ is antimodular, for any embedding of $\mathbb F_2$ into $G$ (since any such embedding is necessarily dense).

The first examples of antimodular actions were found by G. Hjorth \cite{Hj03} who proved that the Bernoulli shift action $\mathbb F_2\curvearrowright \{0,1\}^{\mathbb F_2}$, and its restriction to any invariant set of full product measure, is antimodular. This allowed him to deduce the existence of more than two treeable equivalence relations up to Borel reducibility.
 The result of \cite{Hj03} was later generalized to Borel actions $\Gamma\curvearrowright X$, admitting an invariant probability measure $\mu$ and satisfying certain representation-theoretic conditions. 
Thus, it was shown in \cite{Ke05} that if the Koopman representation $\kappa_0$ of $\Gamma$ on $L^2_0(X)$ is weakly contained in the left regular representation, then $\Gamma\curvearrowright X$ is antimodular. The same conclusion was then shown in \cite{ET07} under the more general assumption that $\kappa_0$ is not amenable in the sense of \cite{Be90}, i.e., that $\kappa_0\otimes\bar{\kappa}_0$ has spectral gap.
On the other hand, if the p.m.p. action $\Gamma\curvearrowright (X,\mu)$ is compact, then $\kappa_0\otimes\bar{\kappa}_0$ admits an abundance of invariant vectors, and thus one cannot apply  the results of \cite{Hj03,Ke05,ET07} to deduce antimodularity. 
As such, Corollary \ref{main} provides the first class of compact actions which are antimodular. 

Although the results of \cite{Hj03,Ke05,ET07} and Corollary \ref{main} cover disjoint classes of actions, as a byproduct of the proof of Theorem \ref{maintech}, 
we are able to generalise the former results as follows. 

\begin{mcor}\label{ET}
Let $\Gamma\curvearrowright (X,\mu)$ be a p.m.p. action admitting a factor action $\Gamma\curvearrowright (Y,\nu)$ such that the Koopman representation $\kappa_0$ of $\Gamma$ on $L^2_0(Y)$ is not amenable. 
Then the restriction of the action $\Gamma\curvearrowright (X,\mu)$ to any measurable $\Gamma$-invariant set $X_0\subset X$ with $\mu(X_0)=1$ is antimodular.
\end{mcor}

In the context of Corollary \ref{ET}, assume moreover that we can decompose $(X,\mu)=(Y\times Z,\nu\times\rho)$, for some probability space $(Z,\rho)$, such that the action $\Gamma\curvearrowright (X,\mu)$ is given by $g\cdot (y,z)=(g\cdot y,c(g,y)\cdot z)$, where $c:\Gamma\times Y\rightarrow$ Aut$(Z,\rho)$ is a cocycle.  (If $\Gamma\curvearrowright (X,\mu)$ is ergodic, then such a decomposition exists by Rokhlin's skew product theorem, see \cite[Theorem 3.18]{Gl03}). In \cite[Theorem 1.2]{ET07} the conclusion of Corollary \ref{ET} is proven under the assumptions that the image of $c$ is contained in a countable subgroup of Aut$(Z,\rho)$. Corollary \ref{ET} removes the countability assumption on $c$.

Finally, we note that Theorem \ref{maintech} leads to a complete characterisation of antimodular ergodic p.m.p.  actions of property (T) groups. Moreover, we show

\begin{mcor}\label{T}
Let $\Gamma$ be a property (T) countable group. Then for any ergodic p.m.p. action $\Gamma\curvearrowright (X,\mu)$, the following conditions are equivalent: 
\begin{enumerate}
\item $\Gamma\curvearrowright (X,\mu)$ is profinite.
\item $\Gamma\curvearrowright (X,\mu)$ is orbit equivalent to a profinite p.m.p. action $\Delta\curvearrowright (Y,\nu)$ of a countable group $\Delta$.
\item There exists  a countable-to-one Borel homomorphism from $\mathcal R(\Gamma\curvearrowright X)_{|D}$ to $\mathcal R(\Delta\curvearrowright Y)$, where  $D\subset X$ is a measurable set with $\mu(D)=1$ and  $\Delta\curvearrowright Y$ is a modular action of a countable group $\Delta$. 
\end{enumerate}
\end{mcor}

Corollary \ref{T} generalizes \cite[Corollary 1.5]{ET07} which shows that  if $\Gamma$ has property (T) and (3) holds, then $\Gamma\curvearrowright (X,\mu)$ is not weakly mixing.

\subsection{Comments on the proof of Theorem \ref{maintech}}  We end the introduction with an informal outline of the proof of Theorem \ref{maintech}. Assume that  $\alpha:X\rightarrow Y$ is a countable-to-one Borel homomorphism from $\mathcal R(\Gamma\curvearrowright X)$ to the orbit equivalence relation of a modular action $\Delta\curvearrowright Y$ of a countable group $\Delta$. 
By \cite{Ke05} we may assume that $\alpha$ is one-to-one.
Our goal is to prove that $\Gamma\curvearrowright X$ is profinite.

Let $\{\mathcal P_n\}_{n\in\mathbb N}$ be a sequence of countable Borel partitions  of $Y$ which separate points such that each $\mathcal P_n$ is $\Delta$-invariant 
and coarser that $\mathcal P_{n+1}$.
Write $\mathcal P_n=\{Y_{n,k}\}_{k\in\mathbb N}$ and define $X_{n,k}=\alpha^{-1}(Y_{n,k})$.
 Then,  as observed in \cite{ET07,Ts09} (see also \cite{Hj03,Ke05}), the partition  $\{X_{n,k}\}_{k\in\mathbb N}$ is almost invariant under the action of $\Gamma$, as $n\rightarrow\infty$. 
 
 We exploit the almost invariance of these partitions to show that the unit vectors $$\eta_n=\sum_{k\in\mathbb N}\frac{1}{\sqrt{\mu(X_{n,k})}}{\bf 1}_{X_{n,k}\times X_{n,k}}\in L^2(X\times X)$$ are almost invariant under the diagonal action of $\Gamma$ on $X\times X$.
 (We make the convention that $\frac{0}{0}=0$, so that only non-negilgible sets $X_{n,k}$ contribute to the definition of $\eta_n$.)
 More precisely, we show that the vectors $\{\eta_n\}$ witness the fact that $\Gamma\curvearrowright X$ is weakly compact (see Lemma \ref{wcomp}). 

Assume for a moment that there exists an ergodic $\Lambda$-invariant measure on $Y$ such that $\alpha$ is measure preserving. Then, after discarding a null set,  each $\mathcal P_n$ is finite and $\Delta$ acts transitively on $\mathcal P_n$. Thus, after relabelling, we have that $\mathcal P_n=\{Y_{n,k}\}_{k=1}^{|\mathcal P_n|}$ and $m(Y_{n,k})=|\mathcal P_n|^{-1}$, for all $k$. In this case, the vectors $\eta_n$ take the simpler form: $\eta_n=\sqrt{|\mathcal P_n|}{\bf 1}_{\sqcup_{k=1}^{|\mathcal P_n|}(X_{n,k}\times X_{n,k})}.$ 
 However, in general, the partitions $\mathcal P_n$ need not be finite, and thus we have to work with the ``weighted" vectors $\eta_n$ defined above.

Returning to our discussion, we combine the spectral gap and the weak compactness properties of $\Gamma\curvearrowright X$ to deduce that this action is compact (see Lemma \ref{compact}). We may therefore identify $\Gamma\curvearrowright X$ with a left translation action $\Gamma\curvearrowright G/H$, where $G$ is a compact group which contains $\Gamma$ densely, and $H<G$ is a closed subgroup. 

By combining again the spectral gap assumption on the action $\Gamma\curvearrowright G/H$ and the $\Gamma$-almost invariance of $\{\eta_n\}$, we conclude that $\eta_n$ is arbitrarily close 
 to a $G$-invariant vector, as $n\rightarrow\infty$. Given $\varepsilon>0$, it follows that for every large enough $n\geq 1$ we have that $$\int_{G}m_{G/H}(g X_{n,k}\cap X_{n,l})^2\;\text{d}m_G(g)\geq (1-\varepsilon)m_{G/H}(X_{n,k})^{3/2}m_{G/H}(X_{n,l})^{3/2},$$ for a ``large proportion" of pairs $(k,l)\in \mathbb N\times\mathbb N$. 
 Thus, if we denote by $\pi:G\rightarrow G/H$ the quotient map, and let $A_{n,k}=\pi^{-1}(X_{n,k})$, then for every large enough $n\geq 1$ we have that $$\int_{G}m_{G}(g A_{n,k}\cap A_{n,l})^2\;\text{d}m_G(g)\geq (1-\varepsilon)m_{G}(A_{n,k})^{3/2}m_{G}(A_{n,l})^{3/2},$$ for a large proportion of pairs $(k,l)\in \mathbb N\times\mathbb N$. 

 In the second part of the proof, we use this condition, and in novel fashion, arguments from the study of approximate subgroups, to show the existence of $k\in\mathbb N$ such that $A_{n,k}$ is ``close" (relative to its measure) to a coset of an open subgroup $G_n<G$ (see Lemma \ref{approx}). Some additional work then shows that the partitions $\{A_{n,l}\}_{l\in\mathbb N}$ and $\{gG_n\}_{g\in G/G_n}$ of $G$ are close, and that $G_n$ contains $H$. 
 Consequently, the partition $\{X_{n,l}\}_{l\in\mathbb N}$ is close to the $\Gamma$-invariant partition $\{gG_nH\}_{g\in G/G_n}$ of $G/H$.
 Since this holds for all $n$ large enough, we conclude that $\Gamma\curvearrowright X$ is profinite.
 
\subsection{Organization of the paper} Besides the introduction, the paper has three other sections. In Section 2, we establish two lemmas needed for the proof of Theorem \ref{maintech}. Section 3 is devoted to the proof of Theorem \ref{maintech}. Finally, in Section 4 we prove Corollaries \ref{main}, \ref{ET} and \ref{T}.

\subsection{Acknowledgements} I would like to thank Lewis Bowen who brought the question addressed by part (1) of Corollary \ref{main} to my attention after learning it himself from Peter Burton.
I am grateful to Alekos Kechris, Todor Tsankov, and Robin Tucker-Drob for stimulating discussions and helpful comments. In particular, I am grateful to Alekos for raising the question which led to part (2) of Corollary \ref{main}.  This work was initiated during the program ``Quantitative Linear Algebra" at the
Institute for Pure and Applied Mathematics. I would like to thank IPAM for its hospitality.
\section{Ingredients for the proof of Theorem \ref{maintech}}

\subsection{Weak compactness} 
The following lemma is the first ingredient in the proof of Theorem \ref{maintech}. Recall that a p.m.p. action $\Gamma\curvearrowright (X,\mu)$ is {\it weakly compact} in the sense of \cite[Definition 3.1]{OP07} if there exists a sequence of non-negative functions $\eta_n\in L^2(X\times X)$ such that 
\begin{enumerate}
\item $\|{\bf 1}_{A\times X}\eta_n\|_2=\|{\bf 1}_{X\times A}\eta_n\|_2=\sqrt{\mu(A)}$, for every measurable set $A\subset X$ and $n\in\mathbb N$.
\item $\lim\limits_{n\rightarrow\infty}\|{\bf 1}_{A\times (X\setminus A)}\eta_n\|_2=0$, for every measurable set $A\subset X$.
\item $\lim\limits_{n\rightarrow\infty}\|\eta_n\circ(g\times g)-\eta_n\|_2=0$, for every $g\in\Gamma$.
\end{enumerate}

\begin{lemma}\label{wcomp}  Let $\Delta\curvearrowright Y$ be a modular action of a countable group $\Delta$ on a standard Borel space $Y$. Assume that $\mathcal P_n=\{Y_{n,k}\}_{k\in \mathbb N}$ is a $\Delta$-invariant Borel partition of $Y$ such that $\mathcal P_{n+1}$ refines $\mathcal P_n$, for every $n\in\mathbb N$, and $\cup_{n\in\mathbb N}\mathcal P_n$ separates points in $Y$.
 Let $\Gamma\curvearrowright (X,\mu)$ be a p.m.p. action. Assume that $\alpha:X\rightarrow Y$ is a Borel map and $X_0\subset X$ is a measurable set such that $\mu(X_0)=1$, the restriction of $\alpha$ to $X_0$ is 1-1, and $\alpha(\Gamma\cdot x)\subset\Delta\cdot\alpha(x)$, for  every $x\in X_0$.

For $n,k\in\mathbb N$, let $X_{n,k}=\alpha^{-1}(Y_{n,k})$ and define $$\eta_n=\sum_{k\in\mathbb N}\frac{1}{\sqrt{\mu(X_{n,k})}}{\bf 1}_{X_{n,k}\times X_{n,k}}\in L^2(X\times X).$$

Then the following hold:
\begin{enumerate}
\item $\|{\bf 1}_{A\times X}\eta_n\|_2=\|{\bf 1}_{X\times A}\eta_n\|_2=\sqrt{\mu(A)}$, for every measurable set $A\subset X$ and $n\in\mathbb N$.
\item $\lim\limits_{n\rightarrow\infty}\|{\bf 1}_{A\times (X\setminus A)}\eta_n\|_2=0$, for every measurable set $A\subset X$.
\item $\lim\limits_{n\rightarrow\infty}\|\eta_n\circ(g\times g)-\eta_n\|_2=0$, for every $g\in\Gamma$.
\end{enumerate}

In particular, the action $\Gamma\curvearrowright (X,\mu)$ is weakly compact.

\end{lemma}

{\it Proof.} Note first that after replacing $X$ with its $\Gamma$-invariant co-null  measurable subset $\cap_{g\in\Gamma}(g\cdot X_0)$, we may assume that $\alpha$ is 1-1 and $\alpha(\Gamma\cdot x)\subset\Delta\cdot\alpha(x)$, for all $x\in X$.

Part (1) is clear. To prove (2), let $A\subset X$ be a measurable set. For $n\geq 1$, let
$e_n$ be the orthogonal projection from $L^2(X)$ onto the $\|.\|_2$-closure of the linear span of $\{{\bf 1}_{X_{n,k}}\}_{k\in\mathbb N}$. Then $$e_n({\bf 1}_A)=\sum_{k\in\mathbb N}\frac{\mu(A\cap X_{n,k})}{\mu(X_{n,k})}{\bf 1}_{X_{n,k}}$$
and a direct calculation shows that $\|{\bf 1}_{A\times (X\setminus A)}\eta_n\|_2=\|{\bf 1}_A-e_n({\bf 1}_A)\|_2$. If we define the Borel partition $\mathcal Q_n=\{X_{n,k}\}_{k\in\mathbb N}$ of $X$, then $\mathcal Q_{n+1}$ refines $\mathcal Q_n$, for every $n\in\mathbb N$. Moreover, since $\alpha$ is 1-1, we have that $\cup_{n\in\mathbb N}\mathcal Q_n$ separates points in $X$. 
As a consequence of these facts, we get that $\|e_n(F)-F\|_2\rightarrow 0$, for every $F\in L^2(X)$, which implies (2).

To prove (3), we follow closely the proof of \cite[Theorem 3.6]{Hj03} (see also \cite[Section 3]{Ke05} and the proof of \cite[Proposition 2.1]{ET07}). Let $g\in\Gamma$ and fix $\varepsilon>0$.

 Let $v, w:X\rightarrow\Delta$ be Borel maps such that $\alpha(g^{-1}\cdot x)=v(x)\cdot\alpha(x)$ and $\alpha(g\cdot x)=w(x)\cdot\alpha(x)$, for almost every $x\in X$. Then we can find a Borel set $A_n\subset X$, for every $n\in\mathbb N$, such that $v$ and $w$ are constant on $A_n\cap X_{n,k}$, for every $k\in\mathbb N$, and $\lim\limits_{n\rightarrow\infty}\mu(A_n)=1$ (see, e.g., the proof of \cite[Lemma 3.9]{Ke05}). Let $\delta_{n,k}\in \Delta$ such that $v(x)=\delta_{n,k}$, for all $x\in A_n\cap X_{n,k}$. 
Since the partition $\mathcal P_n$ is $\Delta$-invariant, we can  define a map $\pi_n:\mathbb N\rightarrow\mathbb N$ such that $\delta_{n,k}Y_{n,k}=Y_{n,\pi_n(k)}.$ If $x\in A_n\cap X_{n,k}$, then  $\alpha(g^{-1}\cdot x)=\delta_{n,k}\alpha(x)\in \delta_{n,k}Y_{n,k}= Y_{n,\pi_n(k)}$, hence $g^{-1}\cdot x\in X_{n,\pi_n(k)}$. This shows that \begin{equation}
\label{theta}g^{-1}\cdot (A_n\cap X_{n,k})\subset X_{n,\pi_n(k)},\;\;\;\text{for all $n,k\in\mathbb N$.}
\end{equation}
Similarly, we can find a map $\sigma_n:\mathbb N\rightarrow N$ such that $g\cdot(A_n\cap X_{n,l})\subset X_{n,\sigma_n(l)}$, for all $n,l\in\mathbb N$.
By combining this and \eqref{theta}, we get that $g\cdot(A_n\cap g^{-1}\cdot(A_n\cap X_{n,k}))\subset X_{n,\sigma_n(\pi_n(k))}\cap X_{n,k}$, for all $n,k\in\mathbb N$.
Thus, if we define $F_n=\{k\in\mathbb N|\mu(A_n\cap g^{-1}\cdot(A_n\cap X_{n,k}))>0\}$, then $\sigma_n(\pi_n(k))=k$, for all $n\in\mathbb N$ and $k\in F_n$.
In particular, the restriction of $\pi_n$ to $F_n$ is $1$-$1$, for every $n\in\mathbb N$.

For $n\in\mathbb N$, let $\tilde\eta_n={\bf 1}_{(g\cdot A_n\cap A_n)\times (g\cdot A_n\cap A_n)}\eta_n$. Since $\lim\limits_{n\rightarrow\infty}\mu(A_n)=1$, by using (1) we get that $\lim\limits_{n\rightarrow\infty}\|\tilde\eta_n-\eta_n\|_2=0$. Thus, in order to prove (3), it suffices to show that $\lim\limits_{n\rightarrow\infty}\|\eta_n-\tilde\eta_n\circ(g\times g)\|_2=0$. To this end, we define $\tilde X_{n,k}=A_n\cap g^{-1}\cdot(A_n\cap X_{n,k})$. Then we have \begin{equation}\label{tildeeta}\tilde\eta_n\circ(g\times g)=\sum_{k\in\mathbb N}\frac{1}{\sqrt{\mu(X_{n,k})}}{\bf 1}_{\tilde X_{n,k}\times\tilde X_{n,k}}.\end{equation}
Now, if $A,B,C\subset X$ are non-negligible measurable sets such that $C\subset A$ and $\mu(C)\leq\mu(B)$, then \begin{equation}\label{comp}\begin{split}&\|\frac{1}{\sqrt{\mu(A)}}{\bf 1}_{A\times A}-\frac{1}{\sqrt{\mu(B)}}{\bf 1}_{C\times C}\|_2^2\\&=\frac{\mu(A)^2}{\mu(A)}-2\frac{\mu(C)^2}{\sqrt{\mu(A)\mu(B)}}+\frac{\mu(C)^2}{\mu(B)}\\&=\frac{\mu(A)^2-\mu(C)^2}{\mu(A)}+\frac{\mu(C)^2}{\mu(A)\mu(B)}\big(\sqrt{\mu(A)}-\sqrt{\mu(B)}\big)^2\\&\leq 2(\mu(A)-\mu(C))+\big(\sqrt{\mu(A)}-\sqrt{\mu(B)}\big)^2\\&\leq 2(\mu(A)-\mu(C))+|\mu(A)-\mu(B)|\\&\leq 2(\mu(A)-\mu(C))+(\mu(A)-\mu(C))+(\mu(B)-\mu(C))\\&=3\mu(A)+\mu(B)-4\mu(C).\end{split}\end{equation}
Since $\{X_{n,k}\}_{k\in\mathbb N}$ is a partition of $X$,  $\tilde X_{n,k}\subset X_{n,\pi_n(k)}$ by \eqref{theta}, and $\mu(\tilde X_{n,k})\leq\mu(X_{n,k})$, for all $n,k\in\mathbb N$, by using \eqref{tildeeta} and \eqref{comp} we deduce that \begin{align*}\begin{split}
&\|\eta_n-\tilde\eta_n\circ(\theta\times\theta)\|_2^2\\&=\sum_{k\in\mathbb N}\|\frac{1}{\sqrt{\mu(X_{n,\pi_n(k)})}}{\bf 1}_{X_{n,\pi_n(k)}\times X_{n,\pi_n(k)}}-\frac{1}{\sqrt{\mu(X_{n,k)}}}{\bf 1}_{\tilde X_{n,k}\times\tilde X_{n,k}}\|_2^2+\sum_{l\in\mathbb N\setminus\pi_n(\mathbb N)}\|\frac{1}{\sqrt{\mu(X_{n,l})}}{\bf 1}_{X_{n,l}\times X_{n,l}}\|_2^2\\ &\leq \sum_{k\in\mathbb N}(3\mu(X_{n,\pi_n(k)})+\mu(X_{n,k})-4\mu(\tilde X_{n,k}))+\sum_{l\in\mathbb N\setminus\pi_n(\mathbb N)}\mu(X_{n,l})\\&\leq 4(\mu(X)-\mu(\cup_{k\in\mathbb N}\mu(\tilde X_{n,k}))\\&=4\mu(X\setminus (A_n\cap g^{-1}\cdot A_n)\end{split}
\end{align*}
Since $\lim\limits_{n\rightarrow\infty}\mu(A_n)=1$, it follows that $\lim\limits_{n\rightarrow\infty}\|\eta_n-\tilde\eta_n\circ(g\times g)\|_2=0$, which finishes the proof.
\hfill$\blacksquare$

\subsection{Approximate subgroups} The following result is the second ingredient needed in the proof of Theorem \ref{maintech}.

\begin{lemma}\label{approx} Let $G$ be a compact  group and denote by $m$ its Haar measure. Let $A, B\subset G$ be non-null measurable sets. Let $\varepsilon\in (0,10^{-5})$ and assume that $$\int_{G}m(gA\cap B)^2\;\text{d}m(g)\geq (1-\varepsilon)m(A)^{3/2}m(B)^{3/2}.$$

Then there exist an open subgroup $G_0<G$ and $g\in G$ such that $m(A\triangle gG_0)<100\sqrt{\varepsilon}\;m(A)$.
\end{lemma}

Lemma \ref{approx} is likely known to the experts but for lack of a reference we provide a complete proof. 
 The proof  of Lemma \ref{approx} is inspired by the proof of \cite[Proposition 5.1]{Ta13}, and by various arguments used in the study of approximate subgroups (see \cite[Chapter 4]{Ta15}). 
Before proving Lemma \ref{approx}, let us record three useful identities that we will use repeatedly. Let $G$ be a compact group, denote by $m$ its Haar measure, and let $A,B\subset G$ be two measurable subsets. Then we have
\begin{equation}\label{1}
\begin{split}
\int_Gm(gA\cap B)^2\;\text{d}m(g)&=\int_{G\times G\times G}{\bf 1}_{gA}(x){\bf 1}_{B}(x){\bf 1}_{gA}(y){\bf 1}_{B}(y)\;\text{d}m^3(g,x,y)\\&=\int_{B\times B}m(xA^{-1}\cap yA^{-1})\;\text{d}m^2(x,y)
 \\&=\int_{B\times B}m(Ax^{-1}\cap Ay^{-1})\;\text{d}m^2(x,y).
\end{split}
\end{equation}

\begin{equation}\label{2}
\begin{split}
\int_Gm(gA\cap B)^2\;\text{d}m(g)&=\int_{G\times G\times G}{\bf 1}_{gA}(x){\bf 1}_{B}(x){\bf 1}_{gA}(y){\bf 1}_{B}(y)\;\text{d}m^3(g,x,y)\\&=\int_{G\times G\times G}{\bf 1}_A(x){\bf 1}_B(gx){\bf 1}_A(g^{-1}y){\bf 1}_B(y)\;\text{d}m^3(g,x,y)\\&=\int_{A\times B}m(yA^{-1}\cap Bx^{-1})\;\text{d}m^2(x,y).\end{split}
\end{equation}

\begin{equation}\label{3}
\begin{split}\int_{A}m(Bx^{-1}\cap B)\;\text{d}m(x)&=\int_{G\times G}{\bf 1}_A(x){\bf 1}_{Bx^{-1}}(y){\bf 1}_{B}(y)\;\text{d}m^2(x,y)\\&=\int_{B}m(A\cap y^{-1}B)\;\text{d}m(y).
\end{split}
\end{equation}

Note also that \eqref{2} implies that 
\begin{equation}\label{4}\int_Gm(gA\cap B)^2\;\text{d}m(g)\leq m(A)m(B)\min\{m(A),m(B)\}\leq m(A)^{3/2}m(B)^{3/2}.
\end{equation}

We first establish a version of Lemma \ref{approx} when the sets in question have equal measures.

\begin{lemma}\label{approx2} Let $G$ be a compact  group and denote by $m$ its Haar measure. Let $A, B\subset G$ be measurable sets such that $\mu:=m(A)=m(B)>0$. Let $\varepsilon\in (0,10^{-4})$ and assume that $$\int_{G}m(gA\cap B)^2\;\text{d}m(g)\geq (1-\varepsilon)\mu^3.$$

Then there exist an open subgroup $G_0<G$ and $g\in G$ such that $m(A\triangle gG_0)<40\sqrt{\varepsilon}\mu$.
\end{lemma}

{\it Proof.} First, note that by \eqref{1} we have that $$(1-\varepsilon)\mu^3\geq \int_Gm(gA\cap B)^2\;\text{d}m(g)=\int_{B\times B}m(Ax^{-1}\cap Ay^{-1})\;\text{d}m^2(x,y)$$
and thus we can find $x_0\in B$ such that $\int_Bm(Ax^{-1}\cap Ax_0^{-1})\;\text{d}m(x)\geq (1-\varepsilon)\mu^2.$

Denote $B_0=\{x\in B\;|\;m(Ax^{-1}\cap Ax_0^{-1})\geq (1-\sqrt{\varepsilon})\mu\}$. Since $m(Ax^{-1}\cap Ax_0^{-1})\leq\mu$, for all $x\in G$, we get that $(1-\varepsilon)\mu^2\leq \int_Bm(Ax^{-1}\cap Ax_0^{-1})\;\text{d}m(x)\leq \mu m(B_0)+(1-\sqrt{\varepsilon})\mu(\mu-m(B_0)),$ which implies that $m(B_0)\geq (1-\sqrt{\varepsilon})\mu$. 
Thus, if we define $A_0=x_0^{-1}B_0$, then $m(A_0)=m(B_0)\in [(1-\sqrt{\varepsilon})\mu,\mu]$ and $m(Ax^{-1}\cap A)\geq (1-\sqrt{\varepsilon})\mu$, for all $x\in A_0$. By using \eqref{3}, this further implies that $$(1-\sqrt{\varepsilon})^2\mu^2\leq \int_{A_0}m(Ax^{-1}\cap A)\;\text{d}m(x)=\int_{A}m(A_0\cap h^{-1}A)\;\text{d}m(h).$$

Therefore, we can find $h\in G$ such that $A_1=A_0\cap h^{-1}A$ satisfies $m(A_1)\geq (1-\sqrt{\varepsilon})^2\mu> (1-2\sqrt{\varepsilon})\mu$.
Define $\varphi:G\rightarrow [0,\infty)$ by letting $\varphi(x)=m(A_1x^{-1}\cap A_1)$, for every $x\in G$.

{\bf Claim.} $G_0=A_1^{-1}A_1$ is an open subgroup of $G$. 

{\it Proof of the claim.} 
As $A_1\subset  h^{-1}A$, we get $m(h^{-1}A\setminus A_1)=m(h^{-1}A)-m(A_1)=\mu-m(A_1)\leq 2\sqrt{\varepsilon}\mu$. If $x\in A_1$, then $x\in A_0$, hence  $m(Ax^{-1}\cap A)\geq (1-\sqrt{\varepsilon})\mu$, and therefore 
\begin{align*}\varphi(x)&\geq m(h^{-1}Ax^{-1}\cap h^{-1}A)-m(h^{-1}Ax^{-1}\setminus A_1x^{-1})-m(h^{-1}A\setminus A_1)\\&=m(Ax^{-1}\cap A)-2m(h^{-1}A\setminus A_1)\\&\geq (1-\sqrt{\varepsilon})\mu-4\sqrt{\varepsilon}\mu=(1-5\sqrt{\varepsilon})\mu.
\end{align*}
Note that $\varphi(x_1^{-1}x_2)\geq\varphi(x_1)+\varphi(x_2)-m(A_1)$, for all $x_1,x_2\in G$. Applying this inequality three times implies that $\varphi(x_1^{-1}x_2x_3^{-1}x_4)\geq\varphi(x_1)+\varphi(x_2)+\varphi(x_3)+\varphi(x_4)-3m(A_1)$, for every $x_1,x_2,x_3,x_4\in G$. 
Since $m(A_1)\leq m(A_0)\leq\mu$ and $\varphi(x)\geq (1-5\sqrt{\varepsilon})\mu$, for every $x\in A_1$, we get that  \begin{equation}\label{gap}\text{$\varphi(y)\geq 4(1-5\sqrt{\varepsilon})\mu-3\mu=(1-20\sqrt{\varepsilon})\mu$, for every $y\in A_1^{-1}A_1A_1^{-1}A_1=G_0 G_0$}.\end{equation}
Since $\varepsilon<10^{-4}$, we deduce that for every $y\in G_0 G_0$ we have $\varphi(y)>0$ and hence $y\in A_1^{-1}A_1=G_0$. Consequently, we get that $G_0G_0\subset G_0$. Since  $G_0\subset G$ is a measurable non-negligible symmetric subset containing the identity, the claim follows. \hfill$\square$

Next, by combining \eqref{3} and \eqref{gap}, we get that \begin{align*}(1-20\sqrt{\varepsilon})\mu\;m(G_0)\leq\int_{G_0}\varphi(y)\;\text{d}m(y)=\int_{A_1}m(G_0\cap z^{-1}A_1)\;\text{d}m(z)\leq m(A_1)^2\leq\mu^2.
\end{align*}
Since $\varepsilon<10^{-4}$, we therefore get that $m(G_0)\leq (1-20\sqrt{\varepsilon})^{-1}\mu<(1+30\sqrt{\varepsilon})\mu$. 

Finally, let $k\in A_1$. Recalling that $G_0=A_1^{-1}A_1$, we get that $A_1\subset kG_0$. Since $A_1\subset h^{-1}A$, we have that $hA_1\subset A\cap hkG_0$, and hence $m(A\cap hkG_0)\geq m(A_1)\geq (1-2\sqrt{\varepsilon})\mu.$ Thus, we have \begin{align*}m(A\triangle hkG_0)&=m(A)+m(hkG_0)-2m(A\cap hkG_0)\\&\leq \mu+(1+30\sqrt{\varepsilon})\mu-2(1-2\sqrt{\varepsilon})\mu< 40\sqrt{\varepsilon}\mu,\end{align*}which finishes the proof.
\hfill$\blacksquare$

\subsection*{Proof of Lemma \ref{approx}} By combining the hypothesis with \eqref{4} we get that $\min\{m(A),m(B)\}\geq (1-\varepsilon)m(A)^{1/2}m(B)^{1/2}$, hence $m(A)/m(B)\in [(1-\varepsilon)^2, (1-\varepsilon)^{-2}]$. Let $\mu=\max\{m(A),m(B)\}$. Let $A_0\supset A, B_0\supset B$ be measurable sets such that  $m(A_0)=m(B_0)=\mu$. Then $$\int_{G}m(gA_0\cap B_0)^2\;\text{d}m(g) \geq (1-\varepsilon)m(A)^{3/2}m(B)^{3/2}\geq (1-\varepsilon)^4\mu^3\geq (1-4\varepsilon)\mu^3.$$

Since $4\varepsilon<10^{-4}$, by applying Lemma \ref{approx} we can find an open subgroup $G_0<G$ such that $m(A_0\triangle gG_0)<40\sqrt{4\varepsilon}\mu\leq 80\sqrt{\varepsilon}(1-\varepsilon)^{-2}m(A)$. Since $m(A\triangle A_0)=\mu-m(A)\leq ((1-\varepsilon)^{-2}-1)m(A)$, we altogether have that $$m(A\triangle gG_0)<(80\sqrt{\varepsilon}(1-\varepsilon)^{-2}+(1-\varepsilon)^{-2}-1)m(A)<100\sqrt{\varepsilon}\;m(A),$$ which proves the conclusion.
\hfill$\blacksquare$

\section{Proof of Theorem \ref{maintech}}\label{section3}

This section is devoted to the proof of Theorem \ref{maintech}. 
Towards proving Theorem \ref{maintech}, we first show that any weakly compact action satisfying the spectral gap condition from its hypothesis is compact.

\begin{lemma}\label{compact} 
Let $\Gamma\curvearrowright^{\sigma} (X,\mu)$ be an ergodic p.m.p. action which is weakly compact.
Let $\kappa$ be the Koopman representation of $\Gamma$ on $L^2(X\times X)$ associated to the diagonal action $\Gamma\curvearrowright (X\times X,\mu\times\mu)$. Denote by $L^2(X\times X)^{\Gamma}\subset L^2(X\times X)$ the subspace of $\kappa(\Gamma)$-invariant vectors.
Suppose that the restriction of $\kappa$ to $L^2(X\times X)\ominus L^2(X\times X)^{\Gamma}$ has spectral gap. Then $\sigma$ is compact.
\end{lemma}

{\it Proof.}
Assume that $\sigma$ is not compact. Let $\Gamma\curvearrowright^{\tau} (Y,\nu)$ be the maximal compact factor of $\sigma$ with the factor map $\pi:(X,\mu)\rightarrow (Y,\nu)$.  
 Then  $L^2(X\times X)^{\Gamma}\subset L^2(Y\times Y)$. 
Consider the embedding $L^{\infty}(Y)\subset L^{\infty}(X)$ given by $f\mapsto f\circ\pi$
and let $E:L^{\infty}(X)\rightarrow L^{\infty}(Y)$ be the conditional expectation.

Since $\Gamma\curvearrowright (X,\mu)$ is weakly compact, we can find a sequence of unit vectors $\eta_n\in L^2(X\times X)$  such that 
$\|{\bf 1}_{A\times (X\setminus A)}\eta_n\|_2\rightarrow 0$ and $\|\eta_n\circ(g\times g)-\eta_n\|_2\rightarrow 0$, for every measurable set $A\subset X$ and $g\in\Gamma$. Since the restriction of $\kappa$ to $L^2(X\times X)\ominus L^2(X\times X)^{\Gamma}$ has spectral gap, we can find a sequence of unit vectors $\xi_n\in L^2(X\times X)^{\Gamma}\subset L^2(Y\times Y)$ such that $\|\eta_n-\xi_n\|_2\rightarrow 0$. Then we have $\|{\bf 1}_{A\times (X\times A)}\xi_n\|_2\rightarrow 0$, for every measurable set $A\subset X$. 

Since $\sigma$ is ergodic and not compact, while $\tau$ is compact, by Rokhlin's skew product theorem (see \cite[Theorem 3.18]{Gl03}) we can decompose $(X,\mu)=(Y,\nu)\times (Z,\rho)$, for some non-trivial probability space $(Z,\rho)$. 
 Let $B\subset Z$ be a measurable set with $0<\rho(B)<1$ and put $A=Y\times B\subset X$. Then $E({\bf 1}_A)=\rho(B){\bf 1}_X$ and $E({\bf 1}_{X\setminus A})=(1-\rho(B)){\bf 1}_X$. 
 Since $\xi_n\in L^2(Y\times Y)$, for every $n$ we have that
 $$\|{\bf 1}_{A\times (X\setminus A)}\xi_n\|_2=\|(E({\bf 1}_A)\otimes E({\bf 1}_{X\setminus A}))\xi_n\|_2=\rho(B)(1-\rho(B))>0,$$
which gives a contradiction.
\hfill$\blacksquare$

{\bf Proof of Theorem \ref{maintech}.} Let $\Gamma\curvearrowright^{\sigma} (X,\mu)$ be an ergodic p.m.p. action satisfying (1) and (2).
Then (2) gives a modular action $\Delta\curvearrowright Y$, a measurable co-null set $D\subset X$, and
a countable-to-1 Borel map $\alpha:D\rightarrow Y$ such that $\alpha(\Gamma\cdot x)\subset\Delta\cdot\alpha(x)$, for all $x\in D$. By \cite[Fact 1.2]{Ke05}, we may assume that $\alpha$ is 1-1. 
 Lemma \ref{wcomp} implies that $\sigma$ is weakly compact. By using (1) and Lemma \ref{compact} we deduce that $\sigma$ is in fact compact. Since $\sigma$ is ergodic, we may assume that $X=G/H$ and $\sigma$ is the left translation action $\Gamma\curvearrowright(G/H,\mu)$, where $G$ is a compact group into which $\Gamma$ embeds densely, $H<G$ is a closed subgroup, and  $\mu=m_{G/H}$. Our goal is to prove that $\sigma$ is profinite.

To this end, it suffices to prove that for any measurable subsets $C_1,...,C_p\subset G/H$ and $\varepsilon>0$, the following holds: we can find a $\Gamma$-invariant measurable partition $\{B_i\}_{i\in I}$ of $G/H$ such that for every $1\leq q\leq p$, there is $I_q\subset I$ satisfying  $\mu(C_q\triangle(\cup_{i\in I_q}B_i))<\varepsilon$.
For the rest of the proof, we fix measurable subsets $C_1,...,C_p\subset G/H$ and $\varepsilon\in (0,1)$. Denote $\delta=\varepsilon^4/10^{12}$.

Let $\mathcal P_n=\{Y_{n,k}\}_{k\in \mathbb N}$ be a $\Delta$-invariant Borel partition of $Y$ such that $\mathcal P_{n+1}$ refines $\mathcal P_n$, for every $n\in\mathbb N$, and $\cup_{n\in\mathbb N}\mathcal P_n$ separates points in $Y$. For $n,k\in\mathbb N$, let $X_{n,k}=\alpha^{-1}(Y_{n,k})\subset G/H$ .
Then $\mathcal Q_n=\{X_{n,k}\}_{k\in\mathbb N}$ is a partition of $D$ such that $\mathcal Q_{n+1}$ refines $\mathcal Q_n$, for every $n\in\mathbb N$. Since $\alpha$ is 1-1, we get that $\cup_{n\in\mathbb N}\mathcal Q_n$ separates points in $D$.

For $n\in\mathbb N$, we define \begin{equation}\label{etan}\eta_n=\sum_{k\in\mathbb N}\frac{1}{\sqrt{\mu(X_{n,k})}}{\bf 1}_{X_{n,k}\times X_{n,k}}\in L^2(G/H\times G/H).\end{equation}

By applying Lemma \ref{wcomp} we derive that $\|\eta_n\circ(g\times g)-\eta_n\|_2\rightarrow 0$, for every $g\in \Gamma$. 
Since the representation of $\Gamma$ on $L^2(G/H\times G/H)\ominus L^2(G/H\times G/H)^{\Gamma}$ has spectral gap by (1), we get $\sup_{g\in\Gamma}\|\eta_n\circ(g\times g)-\eta_n\|_2\rightarrow 0$, as $n\rightarrow\infty$. Since $\|\eta_n\|_2=1$, for all $n\in\mathbb N$, we can find $n\in\mathbb N$ such that 
$\langle\eta_n\circ(g\times g),\eta_n\rangle>1-\delta/2$, for all $g\in\Gamma$. Since $\Gamma<G$ is dense, we deduce that
\begin{equation}\label{uniform}\text{$\langle\eta_n\circ(g\times g),\eta_n\rangle > 1-\delta$, for all $g\in G$}.\end{equation}

 Moreover, since $\cup_{n\in\mathbb N}\mathcal Q_n$ separates points in $D$ and $D\subset G/H$ is co-null, after possibly taking a larger $n$,  we may assume that there exist finite sets $K_q\subset\mathbb N$, $1\leq q\leq p$, such that \begin{equation}\label{A}\text{$\mu(C_q\triangle (\cup_{l\in K_q}X_{n,l}))<\varepsilon/3$, for all $1\leq q\leq p$}.\end{equation}
 For simplicity, since $n$ is fixed for the rest of the proof, we put $X_k=X_{n,k}$ and $\mu_k=\mu(X_{k})\in [0,1]$, for every $k\in\mathbb N$.
Using the definition \eqref{etan} of $\eta_n$, \eqref{uniform} rewrites as
$$\text{$\sum_{k,l\in\mathbb N}\frac{1}{\sqrt{\mu_k\mu_l}}\mu(gX_k\cap X_l)^2> 1-\delta$, for all $g\in G$.}$$
Let $m$ be the Haar measure of $G$. Let $\pi:G\rightarrow G/H$ be the quotient map given by $\pi(x)=xH$. For $k\in\mathbb N$, denote $A_k=\pi^{-1}(X_k)\subset G$. Then $m(A_k)=\mu_k$, for all $k\in\mathbb N$, and the last inequality rewrites as $$\text{$\sum_{k,l\in\mathbb N}\frac{1}{\sqrt{\mu_k\mu_l}}m(gA_k\cap A_l)^2> 1-\delta$, for all $g\in G$.}$$
By integrating over $g\in G$, we derive that $\displaystyle{\sum_{k,l\in\mathbb N}\frac{1}{\sqrt{\mu_k\mu_l}}\int_Gm(gA_{k}\cap A_{l})^2\;\text{d}m(g)> 1-\delta}$. Since $\sum_{k\in\mathbb N}\mu_k=1$, we deduce the existence of $k\in\mathbb N$ such that  \begin{equation}\label{uniform4}\sum_{l\in\mathbb N}\frac{1}{\sqrt{\mu_l}}\int_Gm(gA_{k}\cap A_{l})^2\;\text{d}m(g)> (1-\delta)\mu_k^{3/2}. \end{equation}
Let $S$ be the set of $l\in\mathbb N$ such that $\int_{G}m(gA_{k}\cap A_{l})^2\;\text{d}m(g)> (1-\sqrt{\delta})\mu_k^{3/2}\mu_l^{3/2}$. 
If $b\in S$, then \eqref{4} implies that $\int_Gm(gA_{k}\cap A_{l})^2\;\text{d}m(g)\leq \mu_k^{3/2}\mu_l^{3/2}.$ By combining this inequality with \eqref{uniform4} it follows easily that $\sum_{l\in S}\mu_l\geq 1-\sqrt{\delta}$.
In particular, $S\not=\emptyset$. Since $\sqrt{\delta}=\varepsilon^2/10^6<10^{-5}$, by applying Lemma \ref{approx} to $A_k$ and $A_l$, for some $l\in S$, we deduce the existence of an open subgroup $G_0<G$ and an element $g_k\in G$ such that \begin{equation}\label{uniform5} m(A_{k}\triangle g_kG_0)<100\sqrt[4]{\delta}\mu_k.
\end{equation}
In particular, we get that \begin{equation}\label{G_0}(1-100\sqrt[4]{\delta})\mu_k< m(G_0)< (1+100\sqrt[4]{\delta})\mu_k.\end{equation}

Next, note that if $y\in G_0$, then \eqref{uniform5} implies that $m(A_ky\triangle g_kG_0)=m(A_k\triangle g_kG_0)<100\sqrt[4]{\delta}\mu_k$, and further that $m(A_k\triangle A_ky)<200\sqrt[4]{\delta}\mu_k$. Thus, $m(A_k\cap A_ky)\geq (1-100\sqrt[4]{\delta})\mu_k$, for all $y\in G_0$. By using \eqref{3}, the fact that $G_0$ and $m$ are symmetric,  and \eqref{G_0}, we get that
\begin{equation}\label{uniform6}
\begin{split}\int_{A_k}m(A_k^{-1}x\cap G_0)\;\text{d}m(x)&=\int_{G_0}m(A_k\cap A_ky)\;\text{d}m(y)\\&\geq (1-100\sqrt[4]{\delta})\mu_km(G_0)\\&\geq (1-200\sqrt[4]{\delta})\mu_k^2.
\end{split}
\end{equation}

Now, if $l\in S$ is fixed, then by \eqref{2} we have that \begin{equation}\label{nspe}\begin{split}(1-\sqrt{\delta})\mu_k^{3/2}\mu_l^{3/2}&\leq \int_Gm(gA_k\cap A_l)^2\;\text{d}m(g)\\&=\int_{A_k\times A_l}m(yA_k^{-1}\cap A_lx^{-1})\;\text{d}m^2(x,y)\\&=\int_{A_k\times A_l}m(A_k^{-1}x\cap y^{-1}A_l)\;\text{d}m^2(x,y).\end{split}\end{equation} 
In particular, we get that $(1-\sqrt{\delta})\mu_k^{3/2}\mu_l^{3/2}\leq \mu_k\mu_l\min\{\mu_k,\mu_l\}$, hence $\mu_k/\mu_l\in [(1-\sqrt{\delta})^2,(1-\sqrt{\delta})^{-2}].$
By using \eqref{nspe} again we can find $g_l\in X_{l}$ such that \begin{equation}\label{nspe2}\int_{A_{k}}m(A_{k}^{-1}x\cap g_l^{-1}A_{l})\;\text{d}m(x)\geq (1-\sqrt{\delta})\mu_k^{3/2}\mu_l^{1/2}\geq (1-\sqrt{\delta})^2\mu_k^2.\end{equation}
By integrating the inequality \begin{align*}m(g_l^{-1}A_l\cap G_0)&\geq m((g_l^{-1}A_l\cap A_kx^{-1})\cap (A_kx^{-1}\cap G_0))\\&\geq m(g_l^{-1}A_l\cap A_k^{-1}x)+m(A_k^{-1}x\cap G_0)-\mu_k\end{align*} over $x\in A_k$, and using \eqref{uniform6}, \eqref{nspe2} we get that $$m(A_l\cap g_lG_0)=m(g_l^{-1}A_l\cap G_0)\geq ((1-\sqrt{\delta})^2-200\sqrt[4]{\delta})\mu_k\geq (1-202\sqrt[4]{\delta})\mu_k.$$
Since $\mu_k\geq (1-\sqrt{\delta})^2\mu_l\geq (1-2\sqrt[4]{\delta})\mu_l$, in combination with \eqref{G_0}, we derive that \begin{equation}\begin{split}\label{X_{l}} m(A_{l}\triangle g_lG_0)&=m(A_l)+m(G_0)-2m(A_l\cap g_lG_0)\\&\leq \mu_l+(1+100\sqrt[4]{\delta})\mu_k-2(1-202\sqrt[4]{\delta})\mu_k\\&=\mu_l-(1-504\sqrt[4]{\delta})\mu_k \\&\leq 506\sqrt[4]{\delta}\mu_l,\;\;\;\text{for all $l\in S$}.
\end{split}\end{equation}

We claim that $H\subset G_0$. To see this, let $h\in H$. Since $A_k=\pi^{-1}(X_k)$ we have that $A_kh=A_k$, and by using \eqref{uniform5}, \eqref{G_0} and the fact that $\delta<1/10^{12}$, we get that $$m(G_0h\triangle G_0)=m(g_kG_0h\triangle g_kG_0)\leq 2m(A_k\triangle g_kG_0)\leq 200\sqrt[4]{\delta}\mu_k\leq 200\frac{\sqrt[4]{\delta}}{1-\sqrt[4]{\delta}}m(G_0)<m(G_0).$$
As a consequence we have that $G_0h\cap G_0\not=\emptyset$, hence $h\in G_0$. This shows that indeed $H\subset G_0$.

Thus, if we denote $B=G_0H\subset G/H$ and recall that $A_l=\pi^{-1}(X_l)$, then \eqref{X_{l}} shows that
\begin{equation}\text{$\mu(X_l\triangle g_lB)\leq 506\sqrt[4]{\delta}\mu_l$, for all $l\in S$.}
\end{equation}

Finally,  let $1\leq q\leq p$. Since $\sum_{l\not\in S}\mu_l\leq\sqrt{\delta}$, by using \eqref{X_{l}} we derive that \begin{equation}\begin{split}\label{uniform8}
\mu((\cup_{l\in K_q}X_{l})\triangle (\cup_{l\in K_q\cap S}g_lB))&\leq\sum_{l\in K_q\cap S}\mu(X_{l}\triangle g_lB)+\sum_{l\in K_q\setminus S}\mu(X_{l})\\&\leq\sum_{l\in K_q\cap S}506\sqrt[4]{\delta}\mu_l+\sum_{l\notin S}\mu_l\\& \leq 506\sqrt[4]{\delta}+\sqrt{\delta}.\end{split}
\end{equation}
By combining \eqref{A} and \eqref{uniform8}, we get $m(C_q\triangle(\cup_{l\in K_q\cap S}g_lB))<\varepsilon/3+506\sqrt[4]{\delta}+\sqrt{\delta}\leq \varepsilon/3+507\sqrt[4]{\delta}$. Since $\sqrt[4]{\delta}=\varepsilon/1000$, we conclude that  $m(C_q\triangle(\cup_{l\in K_q\cap S}g_lB))<\varepsilon$, for all $1\leq q\leq p$. Since $\{gB\}_{g\in G}$ is a $\Gamma$-invariant measurable partition of $G/H$, this shows that $\sigma$ is profinite.
\hfill$\blacksquare$

\section{Proofs of Corollaries \ref{main}, \ref{ET} and \ref{T}}
\subsection{Proof of Corollary \ref{main}} Denote by $\kappa$ and $\tilde\kappa$ the Koopman representations of $\Gamma$ associated to the left translation actions $\Gamma\curvearrowright (G,m_G)$ and $\Gamma\curvearrowright (G\times G,m_G\times m_G)$. Since $\Gamma\curvearrowright (G,m_G)$ has spectral gap, we can find $F\subset\Gamma$ finite and $C>0$ such that $\sup_{g\in\Gamma}\|\kappa(g)(\eta)-\eta\|_2\leq C\max_{g\in F}\|\kappa(g)(\eta)-\eta\|_2$, for all $\eta\in L^2(G)$.
Let $U$ be the unitary operator on $L^2(G\times G)$ given by $U(f)(x,y)=f(x,x^{-1}y)$. Then we have that $U^*\tilde\kappa(g)U=\kappa(g)\otimes {\bf Id}_{L^2(G)}$, for all $g\in G$.
Thus, $\tilde\kappa$ is unitarily conjugate  to $\oplus_{i=1}^{\infty}\kappa$, and hence we have that \begin{equation}\label{spec}\text{$\sup_{g\in\Gamma}\|\tilde\kappa(g)(\eta)-\eta\|_2\leq C\max_{g\in F}\|\tilde\kappa(g)(\eta)-\eta\|_2$, for all $\eta\in L^2(G\times G)$.}\end{equation}
If $P:L^2(G\times G)\rightarrow L^2(G\times G)^{\Gamma}$ denotes the orthogonal projection onto the subspace of $\tilde\kappa$-invariant vectors, then $\|P(\eta)-\eta\|_2\leq\sup_{g\in\Gamma}\|\tilde\kappa(g)(\eta)-\eta\|_2$, for all $\eta\in L^2(G\times G)$. In combination with \eqref{spec}, we deduce that $\|\eta\|_2\leq C\max_{g\in F}\|\tilde\kappa(g)(\eta)-\eta\|_2$, for all $\eta\in L^2(G\times G)\ominus L^2(G\times G)^{\Gamma}$. Thus, the restriction of $\tilde\kappa$ to $L^2(G\times G)\ominus L^2(G\times G)^{\Gamma}$ has spectral gap. By applying Theorem \ref{maintech} we get that the action $\Gamma\curvearrowright (G,m_G)$ is profinite, and hence $G$ is a profinite group. \hfill$\blacksquare$

\subsection{Proof of Corollary \ref{ET}} Corollary \ref{ET} is a direct consequence of Lemma \ref{wcomp}. 
Assume by contradiction that the restriction of $\Gamma\curvearrowright (X,\mu)$ to some measurable $\Gamma$-invariant co-null set $X_0\subset X$ is not antimodular. By \cite[Fact 1.2]{Ke05} we can find a modular action $\Delta\curvearrowright Y$ and a 1-1 Borel map $\alpha:X_0\rightarrow Y$ such that $\alpha(\Gamma\cdot x)\subset\Delta\cdot\alpha(x)$, for all $x\in X_0$. By  Lemma \ref{comp} we deduce that $\Gamma\curvearrowright (X,\mu)$ is weakly compact. Since  weak compactness passes to factor actions (see Theorem 6.6 (ii) in the arXiv version of \cite{Io08} or \cite[Proposition 4.9]{Bo17} for a proof of this fact) we deduce that $\Gamma\curvearrowright (Y,\nu)$ is  weakly compact. 
Let $\eta_n\in L^2(Y\times Y)$ be a sequence of unit vectors such that $\|{\bf 1}_{A\times (X\setminus A)}\eta_n\|_2\rightarrow 0$, for all $A\subset Y$ measurable, and $\|\eta_n\circ(g\times g)-\eta_n\|_2\rightarrow 0$, for all $g\in\Gamma$.

 Since $\kappa_0$ is not amenable, $\kappa_0\otimes\bar{\kappa}_0$ has spectral gap. In particular, $\kappa_0$ has spectral gap. By combining these facts with the decomposition $$L^2(Y\times Y)=(L^2_0(Y)\otimes L^2_0(Y))\oplus(\mathbb C{\bf 1}_Y\otimes L^2_0(Y))\oplus(L^2_0(Y)\otimes\mathbb C{\bf 1}_Y)\oplus\mathbb C{\bf 1}_{Y\times Y}$$ of $L^2(Y\times Y)$ into $\Gamma$-invariant subspaces, we deduce the existence of a sequence $c_n\in\mathbb C$ such that  $|c_n|=1$, for all $n$, and $\|\eta_n-c_n{\bf 1}_{Y\times Y}\|_2\rightarrow 0$. Let $A\subset Y$ be a measurable set with $0<\nu(A)<1$.  Then $\|{\bf 1}_{A\times (X\setminus A)}\eta_n\|_2\rightarrow \sqrt{\nu(A)(1-\nu(A))}>0$,  which implies a contradiction. \hfill$\blacksquare$

\subsection{Proof of Corollary \ref{T}} 
It is clear that (1) implies (2) and that (2) implies (3). Assume that (3) holds. Since
 $\Gamma$ has property (T), any unitary representation of $\Gamma$ without non-zero invariant vectors, has spectral gap. Thus, the representation of $\Gamma$ on $L^2(X\times X)\ominus L^2(X\times X)^{\Gamma}$ has spectral gap, and Theorem \ref{maintech} implies (1). \hfill$\blacksquare$

\end{document}